\documentclass[11pt]{amsart}
\usepackage{amssymb,amsmath,amsthm}
\usepackage[hmargin=3cm,vmargin=3.5cm]{geometry}

\usepackage{wrapfig}

\usepackage{graphicx}
\usepackage{amscd}
\usepackage[all, knot]{xy}
\usepackage{epsfig}
\usepackage{psfrag}
\usepackage{amssymb,amsmath,amsthm}
\usepackage{amsfonts}
\usepackage{amssymb}

\newtheorem{theorem}{Theorem}[section]

\newtheorem{problem}[theorem]{Problem}

\newtheorem{exercise}[theorem]{Exercise}

\newcommand\Z{{\mathbb{Z}}}
\newcommand\Q{{\mathbb{Q}}}
\newcommand\R{{\mathbb{R}}}

\newcommand\C{{\mathbb{C}}}

\newcommand\gvect{\mathrm{GVect}}
\newcommand\gdim{\mathrm{gdim}}
\def\stpl{\stackrel{\partial}{\longrightarrow}}
\newcommand{\oplusop}[1]{{\mathop{\oplus}\limits_{#1}}}

\def\lra{\longrightarrow}

\newcommand{\spartial}{\stackrel{\partial}{\lra}}

\title{Linearization and categorification}

\author{Mikhail Khovanov}
\date{March 26, 2016}

\begin{document}
\maketitle

\begin{abstract}
We discuss the notion of linearization through examples, which include the Price
map, PageRank, representation theory, the Euler characteristic and quantum link
invariants. We also review categorification, which adds an additional layer
of structure, in the context of the last two examples.

\end{abstract}


\tableofcontents


\section{Linearization}

We are going to start by discussing \emph{linearization}. We loosely define linearization
as mapping instances of discrete, set-theoretical, geometric, topological, and
even real-life structures (the latter cannot be fully described in a purely mathematical language)
to elements of a linear structure. A linear structure has addition in the background, or,
equivalently, has the structure of an abelian group. Often, a vector space over
a field appears as the background structure, allowing, among other benefits, to use
convenient and powerful tools of linear algebra.

To convince the reader of the ubiquity of linearization let us provide some examples.
We give two examples from real life and several from mathematics.
The real-life examples are

\begin{itemize}
\item  The Price map
\item PageRank
\end{itemize}

\subsection{The Price map}

To naively define the Price map let us start with the semiring of non-negative
integers $\Z_+ = \{ 0, 1, 2, \dots \}$.
Consider the set  Items, which is the union of sets of
Goods and Services:
$$  \mathrm{Items}  \ = \ \mathrm{Goods} \ \cup \  \mathrm{Services}. $$
To value goods and services, in the simplest possible case, we set up a map
$$\mathrm{Price} \ : \ \mathrm{Items} \ \longrightarrow \Z_+.$$
Setting price of a sandwich, for instance, to 6, we value the sandwich at \$6 (using
dollars as the currency, for example).

A more refined version of this setup would map an item together with the data of location
and time to an element of $\Z_+$:
$$ \mathrm{Price} \ : \  (\mathrm{item},\mathrm{time},\mathrm{location}) \ \longrightarrow
\Z_+,$$
since the price of an item might depend on where and when it is offered for sale.

Thousands of books have been written on economics and finance, studying the Price map in
gread detail. This one-page example cannot even attempt to start on Economics 101. Instead,
we would like to point out that the Price map as an
example of \emph{linearization}. Various objects varying from simple to incredibly
sophisticated are mapped to elements of $\Z_+$. Applying this map
loses almost all the information about the objects, retaining only a nonnegative integer.

Most characteristics of an item are lost, but the Price map is incredibly convenient.
It also needs a developed framework to function well (Money, Government, Enforcement).

The structure of the Price map is enhanced by enlarging its target to
\begin{itemize}
\item $\Z$ (with sinister consequences),
\item $\Z\{ \frac{1}{N}\} $ (where $N=100$ is a common choice),
\item A more refined approximation to $\Q$.
\end{itemize}
Enlarging the target of the map to $\Z$, that is, allowing negative prices, is a
natural step from both the real-life and mathematical perspectives. This sometimes leads, in
our interactions with the Price map, to unpleasant consequences,
such as being in debt.

Upon a brief reflection, it becomes clear that the structure of $\Z$ as a commutative
ring is not fully necessary for the Price map. We rarely if ever encounter objects
measured in dollars squared $\$^2$, indicating that the multiplication in the target $\Z$
can be avoided. Rather, it seems enough to make $\Z$, the target of the Price map, a free rank
one module over the ring $\Z$. The natural order $<$ on this module
is, of course, of paramount importance.

That the target is a module rather than a ring becomes further obvious when fractions
are introduced, by enlarging $\Z$ to
$$ \Z \Big \{ \frac{1}{N} \Big \} \ = \  \Big \{ \frac{a}{N}, \ a\in \Z  \Big \} .$$
$N=100$ is a common choice (dollars and cents). Notice that the set $\Z\{\frac{1}{N}\}$
is not naturally a ring. We can view it as a subset of the ring
$$\Z \Big[ \frac{1}{N} \Big] \ = \ \Big \{ \frac{a}{N^n}, \ a\in \Z, n \in \Z_+  \Big \} ,$$
but this subset is not closed under multiplication, indicating, again, that
the target of the Price map is a module rather than a ring.  Thus, the enhanced target
$\Z\{\frac{1}{N}\}$ of the Price map is a $\Z$-module together with the order $<$.

Further enhancement would take us from $\Z\{\frac{1}{N}\}$ to $\Q$ (viewed
as a module over itself), or even to $\R$,
with the natural order $<$ extended to $\Q$ and $\R$.
Various approximations to $\Q$ and $\R$, such as the floating point type, are used
as practical implementations of this enhancement.

In all cases, the target carries a natural linear structure
(abelian semigroup $\Z_+$ in the initial example, and abelian groups $\Z$, $\Z\{\frac{1}{N}\}$,
and $\Q$ in the generalizations).

\begin{problem}
Find and implement modifications of the Price map to other, more refined, targets.
\end{problem}

For instance, is there a useful modification such that the analogue of the addition
operation on the target is noncommutative? Is there a modification where the target
is (a free module over) a noncommutative ring?

\vspace{0.02in}

Bitcoin is a recent perplexing example of the Price and Money framework
where the target structure is a subset of $\Z_+$
that monotonically grows with time and was designed to be bounded from above at all times by
approximately 20.3 million~\cite{N}.  Bitcoins
are awarded for solving hard meaningless instances of a  computational problem, increasing in complexity as more bitcoins
are minted, which results in the above upper bound.  This artificial upper bound on the total number of bitcoins
is only partially resolved by introducing fractional bitcoins.

\subsection{PageRank}

Already back in the late 90's the structure of the Internet was incredibly complex.
The PageRank idea~\cite{PSMW}
was to take the Internet - the set of all webpages together with
their content and data - and forget almost everything about it, reducing it
to a directed graph, the web graph.
Vertices of the web graph $\Gamma$ are webpages and there is
an oriented edge from $i$ to $j$ if $i\not= j$ and there is a link from $i$ to $j$.

This first step converts an enormously complicated structure to a simple one,
at least on the theoretical level. The second step is linearization. We form the
web matrix $A$. Its rows and columns are labelled by
vertices, and the entry
$$a_{ji} = \frac{1}{N_i}$$
if there is an edge from $i$ to $j$, where $N_i$
is the number of oriented edges out of $i$. Otherwise $a_{ji} = 0$.

Assuming irreducibility (there is a path from any $i$ to any $j$), the web matrix
has the maximum eigenvalue and unique (up to scaling) eigenvector with this
eigenvalue - the PageRank vector $v$. In greater generality, this is known as the
Perron-Frobenius eigenvalue and eigenvector.

In this linearization we form a real vector space $V$ with the basis labelled by
webpages and convert the web graph $\Gamma$ to a linear transformation on $V$
given by the matrix $A$. Coefficients of the maximum eigenvector $v$ are
positive and rank the webpages by a measure of their popularity.

Just like the Price map, this conversion of the entire Internet to a linear
transformation and its maximum eigenvector loses almost the entire informational
content of the Internet, yet it proved to carry phenomenal value.

\vspace{0.2in}

We conclude our discussion of real-life linearizations with an exercise
and a project for the reader.

\begin{exercise}
Find more examples of real-life linearizations.
\end{exercise}

\begin{problem}
Discover a new real-life linearization and develop it.
\end{problem}

\subsection{Representation theory as linearization}

The notion of a group $G$ acting on a set $X$ is one of the earliest
fundamental concepts we encounter in modern algebra. A very natural
generalization of a (left) group action on a set is that of a monoid
 (a semigroup with the unit element) $G$ acting on a set $X$, via a map
$$  G \ \times \ X \ \lra \ X $$
subject to the associativity and unitality constraints.

The notion of an action of a monoid on a set admits a linearization.
Linearization of a monoid is a ring $R$, linearization of a set on
which a monoid acts is that of a (left) module $M$ over $R$, with
the action being a bilinear map
$$ R \ \times  \ M \ \lra \ M$$
subject to the usual axioms. Study of modules over rings, also
known as representation theory, is of utmost importance in
modern mathematics.

Transformation from $(G,X)$ to $(R,M)$ can be achieved in two
steps. First, linearize monoid $G$ to the semigroup algebra $FG$,
where $F$ is a field and elements of $FG$ are finite linear
combinations of elements of $G$ with coefficients in $F$.
Linearize $X$ to the $F$-vector space $FX$ with basis $X$.
Algebra $FG$ acts on its module $FX$. This action is often
much more interesting than the corresponding action of $G$ on
$X$, as we can see, for example, from the case when $G$ is the
symmetric group $S_n$. Orbits of group-theoretical actions of $S_n$ on
sets are classified by conjugacy classes of subgroups of $S_n$, while
actions of $FS_n$ on $F$-vector spaces (representations of
the symmetric group) constitute
 a beautiful theory with many applications to geometry, topology
 and algebra.

Furthermore, some parts of representation theory, especially representations
of simple Lie algebras, quantum groups, and Hecke algebras, can
be categorified. This note does not discuss such categorifications,
instead restricting to more topological examples that include
the Euler characteristic of a topological space and the Jones
polynomial.

\section{Euler characteristic and homology of topological spaces}

The Euler characteristic can be thought of
 a map from the set of sufficiently nice topological
spaces to the ring  of integers
$$ \chi \ : \ \mathrm{Nice \ topological \ spaces} \ \lra \ \Z $$
At first, it is defined naively for only, say, finite simplicial complexes.
The latter are spaces given by starting with finitely many disjoint
simplices of various
dimensions and identifying their subsimplices via linear maps. If the resulting
simplicial complex $M$ has $|M|_n$ simplices of dimension $n$, the Euler
characteristic of $M$ is
$$\chi(M) \ = \ |M|_0 - |M|_1 + |M|_2 - \dots = \sum_n (-1)^n |M|_n ,$$
the alternating sum of the number of simplices in $M$ of various dimensions.

\begin{theorem} $\chi(M)$ is an invariant of $M$.
\end{theorem}

There are various ways in which $\chi$ is an invariant. The simplest way to
phrase the invariance is as the independence of $\chi(M)$ on the choice of the
simplicial decomposition of $M$, but, in fact, it has a much stronger invariance property, depending only on the homotopy type of $M$.

We can think of the Euler characteristic map $\chi$ as a linearization.
A set of sufficiently well-behaved topological spaces, those that have
a realization as finite simplicial complexes, is mapped to integers.
Almost all topological information about a space is lost under this
map, but it does provide an invariant of the space. Simple operations
of the disjoint union and direct product of spaces correspond to addition
and multiplication on their Euler characteristics:
$$ \chi(M\sqcup N) = \chi(M) + \chi(N), \ \ \chi(M\times N) = \chi(M) \chi(N).$$

The Euler characteristic is a useful but
rather basic invariant of topological spaces. It can be enhanced by
keeping track of the boundaries of simplices. The boundary of an
$n$-simplex is a union of its $n+1$ facets, each an $(n-1)$-simplex.
The boundary operation can be though of as assigning to an $n$-simplex
the set of its facets, subsimplices of codimension $1$.
We linearize the boundary operation by taking a simplex $v$ to the linear
combination of its facets, with $\pm 1$ coefficients, that is, signs
$$ \partial(v) \ = \ \sum \ \pm \mathrm{face}(v).$$
Signs come from keeping track of orientations of $v$ and its facets.

The intuitive observation
$$\emph{The \ boundary \ of \ the \ boundary \ is \  empty}$$
linearizes to
$$ \partial \partial = 0 $$
or $\partial^2=0$. This happens since each $(n-2)$-simplex on the boundary
of $v$ appears twice in the expansion of $\partial^2(v)$, with opposite signs,
leading to $\partial^2(v)=0$ for all $v$, thus $\partial^2=0$.

To make this more formal, let $S_n$ be the set of all $n$-simplices of the
decomposition, and $V_n = \Q\langle S_n \rangle$, a $\Q$-vector space with basis $S_n$.
(Instead of $\Q$ one can use $\Z$, $\R$, or any commutative ring).
Extend $\partial$ in a linear fashion to a map $V_n \lra V_{n-1}$.
The result is a complex
$$ \dots \spartial V_{n+1}\spartial V_n \spartial V_{n-1} \spartial \dots .$$
Since $\partial^2=0$, we can define the $n$-th homology groups of the complex $V$:
$$ H_n(V) \ = \ (\mathrm{ker} \ \partial \ : \ V_n \lra V_{n-1}) \ / \ (\mathrm{im}\ \partial \ : \ V_{n+1} \lra
V_n),$$
by taking the quotient of the subspace of $V_n$ which is the kernel of the boundary map
by the image of $V_{n+1}$ under $\partial$.

This construction allows one to define homology groups $H_n(M)$ of a simplicial
complex $M$ as $H_n(V)$, for the above complex $V$ of vector spaces.
An important result establishes that $H_n(M)$ are isomorphic to homology
groups (singular homology)
defined in a more invariant way, via the complex generated
in degree $n$ by all continuous maps from a fixed $n$-complex to $M$,
and the boundary operator given by essentially the same formula as above.

Singular homology, which we also denote $H_n(M)$, has the benefit of
being defined for any topological space $M$. To $M$ we associate the
total homology groups
$$H(M) \ = \ \oplusop{n\ge 0} H_n(M).$$
$H(M)$ is a graded vector space (or an abelian group, if the ground
ring in $\Z$).

\begin{theorem} The Euler characteristic
$$ \chi(M) \ = \ \sum_{n \ge 0} (-1)^n \dim \ H_n(M).$$
\end{theorem}

Thus, we can recover the Euler characteristic $\chi$ from
a more refined invariant - simplicial or singular homology
groups. The benefits of homology groups over the Euler
characteristic are plentiful; what follows
is a rather incomplete list.

\begin{enumerate}
\item Singular homology groups are defined for all
topological spaces, not only nice ones (which in our
case meant finite simplicial complexes).
\item They carry much higher informational content than
the Euler characteristic.
\item Homology groups are \emph{functorial}. To a continuous
map of topological spaces $f:X\lra Y$
there is associated a homomorphism of group (or vector spaces)
$$f_{\ast} \ : \ H_n(X) \lra H_n(Y)$$
These homomorphisms together form a functor
$$ \mathrm{Top} \stackrel{H}{\lra} \mathrm{GrAb}$$
 from the category
of topological spaces and continuous maps to the category
of graded abelian groups.
\item Any commutative ring can be used for coefficients (the
most common choices are $\Z,\Q,\Z/n\Z,\R,\C$).
\item Beyond functoriality, homology carries additional
structures. One of them (comultiplication)
is easier to understand on the dual object
to $H_{\ast}(M)$, the cohomology groups $H^{\ast}(M)$,
which are naturally a graded super-commutative ring for any
$M$.
\item Homology and cohomology generalize to exraordinary
(co)homology theories, carrying even more information.
\end{enumerate}

Passing from the Euler characteristic to homology of topological
spaces is an example of \emph{categorification}. Homology
functor
$$ \mathrm{Top} \stackrel{H}{\lra} \mathrm{GrAb}$$
lifts the Euler characteristic (linearization) map
$$\mathrm{(Nice)\ Topological \ Spaces} \ \stackrel{\chi}{\lra}
\ \mathrm{Integers}$$
We now provide a very basic categorification dictionary
which helps to explain the lifting of integer-valued invariants
that categorification provides.

\begin{table}
\begin{tabular}{c | c | c
 }
Structure  & Elements and      & Categorification  \\
           & operations        &                   \\
\hline \hline
       &   $n,m$ &  vector spaces $V$, $W$  \\
$\Z_+$ &  $n+m$  &  direct sum $V\oplus W$  \\
       &  $n\cdot m$ & tensor product $V\otimes W$ \\
\hline
       &   $n,m$   &  complexes $V$,$W$ of vector spaces \\
$\Z$   &   $n-m$   &  cone of a map $f:W \lra V$ \\
       &          &       \\
 \hline
   &       & It is an open problem \\
  $\Q$      &    $n/m$    &  to  categorify division \\
        &        &
\end{tabular}
\end{table}

Let us start with the first of the three rows in this
dictionary.
Category $\mathrm{Vect}$ of finite-dimensional vector spaces
over a field categorifies the semiring $\Z_+$ of nonnegative integers. To an object $V$ of $\mathrm{Vect}$ we assign its dimension
$\dim(V) \in \Z_+$. Direct sum and tensor product of vector
spaces decategorify to the addition and multiplication of
nonnegative integers
$$ \dim(V\oplus W) = \dim(V) + \dim(W), \ \
 \dim(V\otimes W) = \dim(V) \otimes \dim(W). $$
 Upon decategorification, all information about the morphisms
 (and most information about the objects) is lost.

 In this model example, decategorification maps objects of the
 additive monoidal category $\mathrm{Vect}$ to elements of
 the semiring $\Z_+$.
Instead of $\mathrm{Vect}$ the category of finitely-generated
free abelian groups can be used as well, with the rank of
the free group taking place of dimension.

The semiring $\Z_+$ naturally sits inside the ring $\Z$.
 To lift subtraction of integers to the categorical level,
 we enlarge the category from vector spaces to that of complexes
 of vector spaces. A possible natural restriction is to require that
 complexes are bounded and finite-dimensional in each degree.
 This can be relaxed to requiring that the total homology groups
 of a complex are finite-dimensional, so that the Euler
 characteristic of an object in this category is well-defined.

At first, the morphisms in this category are just homomorphisms of
complexes (linear maps of vector spaces in each degree that
intertwine differentials in the two complexes). To get a useful
category, morphisms need to be modified -- we mod out by the
ideal of null-homotopic morphisms.
The resulting category $\mathcal{C}$ of complexes modulo chain homotopies
is \emph{triangulated}, and its Grothendieck group is the
corresponding group of a triangulated category. Tensor product
of vector spaces naturally extends to complexes, respecting the
ideal of null-homotopic morphisms, and leads to a tensor structure
on $\mathcal{C}$. The Grothendieck group $K_0(\mathcal{C})$ of $\mathcal{C}$
acquires the structure of a ring, and there is a natural ring isomorphism
$$K(\mathcal{C}) \ \cong \ \Z .$$
This isomorphism is induced by the map that
takes an object $V$ of $\mathcal{C}$ to its Euler characteristic $\chi(V)$.

The analogue of the subtraction operation on integers is the cone
of a map of complexes. Given a map $f: V \lra W$ of complexes,
shift the source complex $V$ one step to the left, form direct
sums of vector spaces $V^{n-1}\oplus W^n$, over all $n$, and define
the differential in the new complex as $-d_V+d_W+f$. The cone
complex $\mathrm{Cone}(f)$ has Euler characteristic the difference of
those for $W$ and $V$:
$$\chi(\mathrm{Cone}(f)) \ = \ \chi(W) - \chi(V).$$

Category $\mathcal{C}$ is one of the most fundamental and useful
monoidal triangulated categories with the Grothendieck group $\Z$.
In the hierarchy of structures we built out of integers the
next object, in complexity, is the ring of rational numbers $\Q$,
which, unlike integers, allows division by a nonzero number.
At this point we are already at a limit of current mathematical
knowledge - it is not known how to categorify division and the
ring of rational numbers $\Q$. We record this as an open problem.

\begin{problem}
Describe a triangulated monoidal category $\mathcal{C}$ with
the Grothendieck ring isomorphic to $\Q$.
\end{problem}

Even the following apparently simpler problem appears to be open (an approach
to the $n=2$ case is being considered in \cite{KT}).

\begin{problem}
Construct a triangulated monoidal category $\mathcal{C}$ with
the Grothendieck ring isomorphic to $\Z\big[ \frac{1}{n}\big]$.
\end{problem}

\section{Jones polynomial and its categorification}

\subsection{Jones polynomial and the Kauffman bracket}
Knots and links are smooth or piecewise-linear embeddings
of a single circle $S^1$ (knots) or a disjoint union of finitely-many
circles (links) into $\R^3$, with the embeddings considered up
to isotopies. They appear toyish at first glance, but in the
past few decades have been related to an amazing plethora
of deep structures in mathematics and mathematical physics.

On such structure is the Jones polynomial \cite{J}. It is an
invariant of links that assigns a Laurent polynomial $J(L)$
in a single variable $q$ to an oriented link in $\R^3$, and
can be thought as a map
$$ J \ : \ \mathrm{Links} \ \lra \ \Z[q,q^{-1}] $$
from the set of oriented links in $\R^3$ to a linear structure -
the ring of Laurent polynomials $\Z[q,q^{-1}]$.

The Jones polynomial is uniquely determined by the conditions:
\begin{itemize}
\item The Jones polynomial of the trivial knot is $q+q^{-1}$,
\item For any three links that differ only in the
neighbourhood of a small ball as depicted below, there is
a linear relation on their Jones polynomials
\begin{equation}\label{skein}
q^2 J \left( \raisebox{-0.3cm}{\psfig{figure=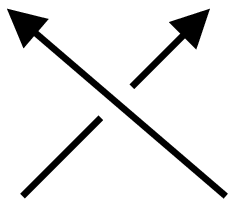,height=0.8cm}}\right) -
q^{-2} J \left( \raisebox{-0.3cm}{\psfig{figure=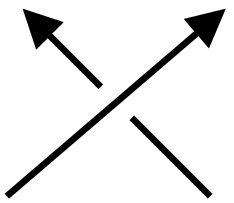,height=0.8cm}}\right)
= (q-q^{-1})J\left( \raisebox{-0.3cm}{\psfig{figure=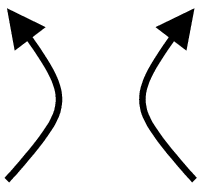,height=0.8cm}}\right).
\end{equation}
\end{itemize}

Shortly after the discovery of the Jones polynomial, Louis Kauffman~\cite{Ka}
found a recursive construction of the polynomial that allows
an elementary proof that the polynomial is well-defined and
provides a wealth of other structural information. His construction
is known as the Kauffman bracket.

To define the Kauffman bracket, take a
generic projection $D$ of an oriented link $L$
onto the plane $\R^2$ (generic in the sense of not having triple
intersection points and tangency points) and temporarily forget
about the orientation of $L$.

If projection $D$ has no crossings, we define its Kauffman
bracket
$$\langle D \rangle = (q+q^{-1})^c,$$
where $c$ is the number of components (circles) in the projection.

\begin{figure}[!htbp]
\includegraphics[width = 0.2\textwidth]{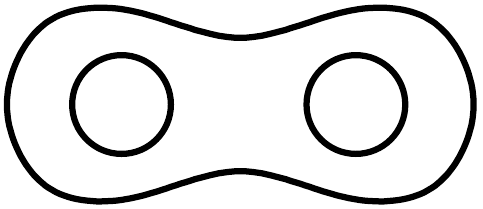}
\raisebox{0.48cm}{$ \quad \longmapsto \quad (q+q^{-1})^{\# \ of \ circles}   $}
\end{figure}

$D$ having crossings is an interesting case.
We pick a crossing of $D$ and define the bracket of $D$
recursively, as  a linear combination of brackets of projections
with one crossing less:
  \begin{figure}[!htbp]
    \includegraphics[width = 0.06\textwidth]{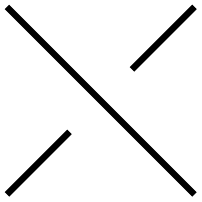}
\raisebox{0.3cm}{$ \ \ = \ \ $}  \includegraphics[width = 0.055\textwidth]{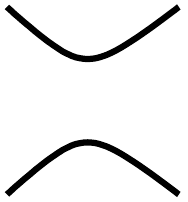}
\raisebox{0.3cm}{$\ \  - \ q^{-1} \ \ $}  \includegraphics[width = 0.065\textwidth]{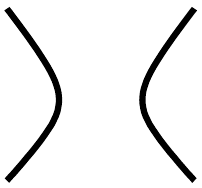}
  \end{figure}

If projection $D$ has $n$ crossings, the full expansion of $\langle D \rangle$ will have
$2^n$ terms, each of the form $ \pm q^{a}(q+q^{-1})^{b}$
for some $a\in \Z$ and $b \in \Z_+.$ This construction uniquely
determines $\langle D \rangle$ for any generic diagram $D$.
To get a link invariant, recall that $L$ is oriented and
define
\begin{equation}\label{jones-kauffman}
J(D) \ = \ (-1)^{x(D)} q^{2x(D)-y(D)}  \langle D \rangle ,
\end{equation}
where $x(D)$ and $y(D)$ is the number of negative
 \raisebox{-0.16cm}{\psfig{figure=m1.eps,height=0.5cm}} and
 positive \raisebox{-0.16cm}{\psfig{figure=m2.eps,height=0.5cm}}
crossings of $D$.


  \begin{theorem} (L.Kauffman) The resulting polynomial $J(D)$ does not depend
  on a choice of projection $D$ of a link $L$ and equals the
  Jones polynomial $J(L)$.
  \end{theorem}

  The theorem claims the invariance of $J(D)$ under the Reidemeister
  moves of link diagrams and has a direct computational proof.
  Kauffman's bracket gives the easiest way to see that the
  Jones polynomial is well-defined.

\subsection{Graded complexes}
We now move on to categorification of the Jones polynomial. It takes
values in $\Z[q,q^{-1}]$, and we start by realizing this ring
as the Grothendieck ring of a suitable category. Recall that
the ring $\Z$ was lifted, at first, to the category of
finite-dimensional vector spaces, and later, to the category of
complexes of vector spaces. Analogous lifting of $\Z[q,q^{-1}]$
is realized via graded vector spaces. Consider the category
$\gvect$ with objects - graded finite-dimensional vector spaces
(say, over $\Q$)
$$ V \  = \ \oplusop{n\in \Z} V_n \ = \ \quad \ldots \oplus V_{-2} \oplus V_{-1} \oplus V_0 \oplus V_1
\oplus V_2 \oplus \dots $$
In particular, only finitely many of $V_n$'s are nonzero.
Morphisms in $\gvect$ are linear maps $V\lra W$ that preserve
the grading, that is, take $V_n$ to $W_n$ for all $n$.
To $V$ assign its graded dimension
$$\gdim(V) \ = \ \sum_{n\in \Z}\ \dim(V_n)\cdot q^n \in \Z[q,q^{-1}].$$
One can "add" and "multiply" graded vector spaces, by forming
their direct sum and tensor product. These operations turn
$\gvect$ into a linear monoidal category. Its Grothendieck ring
is naturally $\Z[q,q^{-1}]$, and variable $q$ becomes a grading
shift upon this lifting.

$\mathrm{gdim}(V),$ for a graded vector space $V$, has non-negative
coefficients. To allow arbitrary integer coefficients we need
to pass to complexes, while maintaining the grading.
A complex of graded vector spaces $V$ is a bigraded
vector space
$$ V \ = \ \oplusop{n,m\in \Z} V^m_n$$
with a differential $\partial: V \lra V,$ $\partial^2=0$, that respects
additional grading, that is, restricts to
$$\partial: V^m_n  \lra    V^m_{n-1}$$
for all $n,m$. We can think of $V$ as the direct sum of complexes of
vector spaces
$$V^m \ \ = \ \ \dots \stpl V^m_{n+1} \ \stpl V^m_n \ \stpl V^m_{n-1} \  \stpl \dots $$
in each $q$-degree $m$.

To a complex of graded vector spaces $V$ we associate homology groups
$$ H_n^m(V) \ = \ (\mathrm{ker} \ \partial \ : \ V^m_n  \longrightarrow    V^m_{n-1} \ ) /
\ (\mathrm{im} \ \partial \ : \ V^m_{n+1} \longrightarrow V^m_n \ ) $$
Homology of $V$ is a bigraded vector space
$$ H(V) \ = \ \bigoplus_{n,m\in \Z} \ H_n^m(V) .$$
As for complexes, we pick a suitable category $\mathcal{C}$
to work with, by
requiring that homomorphisms and homotopies of graded complexes
respect the extra grading $m$, and imposing the finite-dimensionality
condition $\dim(H(V))< \infty$.

For an object $V$ of $\mathcal{C}$ the Euler characteristic of each
 complex  $V^m= \oplusop{n\in \Z}V^m_n$ is an integer
$$ \chi(V^m) \ = \ \sum_n (-1)^n \dim (H_n^m(V)).$$
If $V^m$, and not just its homology, is finite-dimensional,
the Euler characteristic can be computed from the spaces themselves,
$$ \chi(V^m) \ = \ \sum_n (-1)^n \dim (V^m_n).$$

We make these integers into coefficients of a Laurent polynomial
$$\chi(V) \ = \ \sum_m \chi(V^m)\cdot q^m \ = \ \sum_{m,n}
\ (-1)^n \dim (H^m_n(V)) \cdot q^m. $$

When $V$ is finite-dimensional, we also have
$$\chi(V) \ = \ \sum_{n,m} \ (-1)^n \dim(V^m_n) \cdot q^m $$

\subsection{Categorification}

The idea behind categorification of the Jones polynomial $J(L)$
is to look at the formula for the Kauffman bracket $\langle D
\rangle$ of a diagram $D$ and consistenly lift all the
terms there into a complex $C(D)$ of graded vector spaces with
the Euler characteristic $\langle D \rangle$.

Key term in the construction is $q+q^{-1}$, the Kauffman bracket of
a simple circle in the plane, also equal to the Jones polynomial
of the unknot. We lift this polynomial to the graded
vector space
$$A \ = \ \Q \cdot \mathbf{1} \oplus \Q \cdot X $$
with basis elements denoted $\mathbf{1}$ and $X$, with degrees
$$ \deg \ \mathbf{1} \ = \  -1 , \ \  \deg \ X \ = \ 1.$$
 Graded dimension of $A$ is $q+q^{-1}$. The complex associated
 to the simple circle diagram will be just $A$, placed in
 homological degree $0$:
 $$ 0 \lra A \lra 0.$$

 To a crossingless diagram with two circles we assign $A\otimes A$,
  which has a basis
  $ \{ \mathbf{1} \otimes \mathbf{1}, \  \mathbf{1} \otimes X, \
  X \otimes \mathbf{1}, \ X\otimes X \} .$
In general, to a diagram which consists of $k$ disjoint, perhaps
nested, circles we associate $A^{\otimes k}$, which has
the graded dimension $(q+q^{-1})^k $.

We next move on to diagrams with a single crossing. For such
a diagram $D$ its Kauffman bracket $\langle D\rangle$ is
the difference of two terms, one coming from a two-circle
diagram, the other from a one-circle diagram, corresponding to
two different ways to simplify (resolve) $D$ into crossingless diagrams.
One such example is depicted below.

  \begin{figure}[!htbp]
 \includegraphics[width = 0.4\textwidth]{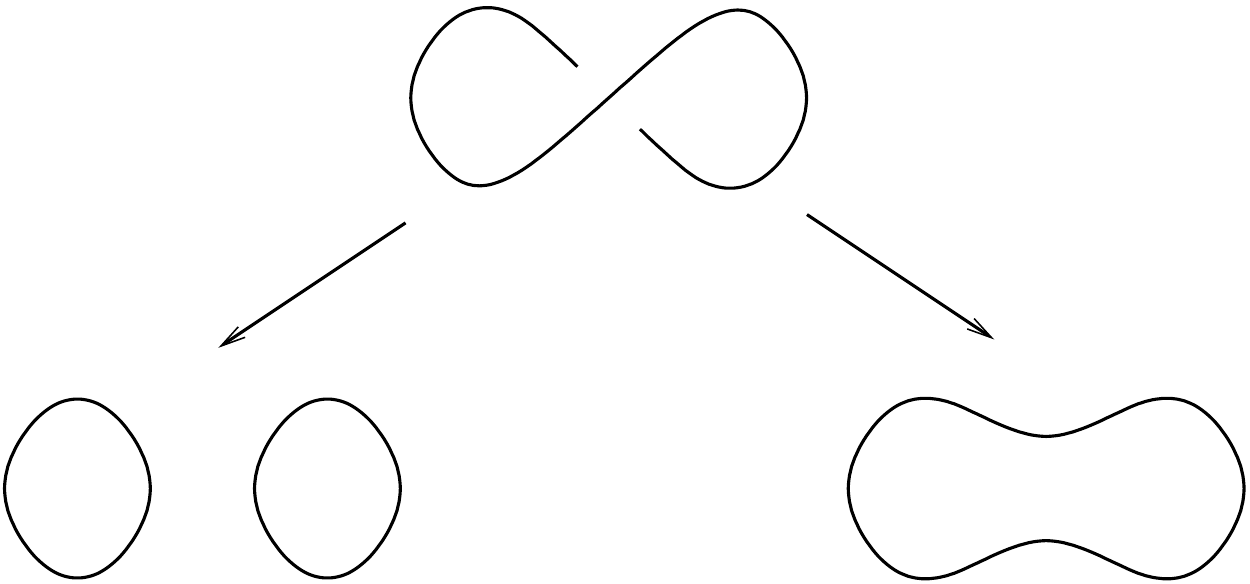}
 \end{figure}

 In the Kauffman formula for $\langle D\rangle$, the single circle
diagrams enters with the coefficient $-q^{-1}$, and the two-circle
one with coefficient $1$. To interpret $\langle D\rangle$ as
the Euler characteristic, we place $A^{\otimes 2}$ in homological
degree $0$, $A$ in homological degree $1$, and look for the
differential to make this into a complex.
 $$ 0 \lra \ A  \otimes A \ \stackrel{M}{\lra} \ A \  \lra \  0 $$
We label the differential $M$, for \emph{multiplication}, since
that's what it looks like. To interpret $q^{-1}$ in the formula,
we need to shift the internal grading of $A$ down by $1$.
Since we want homology to be an invariant of a knot, and not just
its diagram, the homology
of the above complex should be isomorphic to $A$, perhaps up to
an overall grading shift, that we can take care of later.
Therefore, $M$ is surjective and, moreover,
it must preserve the internal degree. The table below
lists the degrees of basis elements of the two spaces, with
$\{-1\}$ denoting the degree shift down by $1$.

 \begin{table}[!htbp]
 \begin{tabular}{c|c|c}
 Degree  &  Basis of $A\otimes A$ & Basis of $A\{-1\}$  \\
 \hline \hline
  $2 $  & $ X \otimes X$   &     \\
 \hline
  $1$   &                & \\
 \hline
 $0$   & $ X\otimes \mathbf{1}, \mathbf{1}\otimes X$ & $ X$  \\
 \hline
  $-1$  &              & \\
 \hline
 $-2$   &  $ \mathbf{1} \otimes \mathbf{1} $ &  $ \mathbf{1}$
 \end{tabular}
 \end{table}

\newpage

$M$ is now constrained to essentially a unique map, given below.
 \begin{table}[!htbp]
 \begin{tabular}{lll}
 $X \otimes X$   &  $\longmapsto$ & $0$ \\
 $X\otimes \mathbf{1}, \  \mathbf{1}\otimes X$ & $\longmapsto$ & $X$ \\
 $\mathbf{1} \otimes \mathbf{1}$ &  $\longmapsto$ & $\mathbf{1}$
 \end{tabular}
 \end{table}

This map makes $A$ into a commutative associative algebra with
the unit element $\mathbf{1}$.

The other case for a single crossing complex is when the $-q^{-1}$
coefficient in the Kauffman formula appears with the $(q+q^{-1})^2$ term. A possible diagram when this happens is depicted below.

 \begin{figure}[!htbp]
\includegraphics[width = 0.4\textwidth]{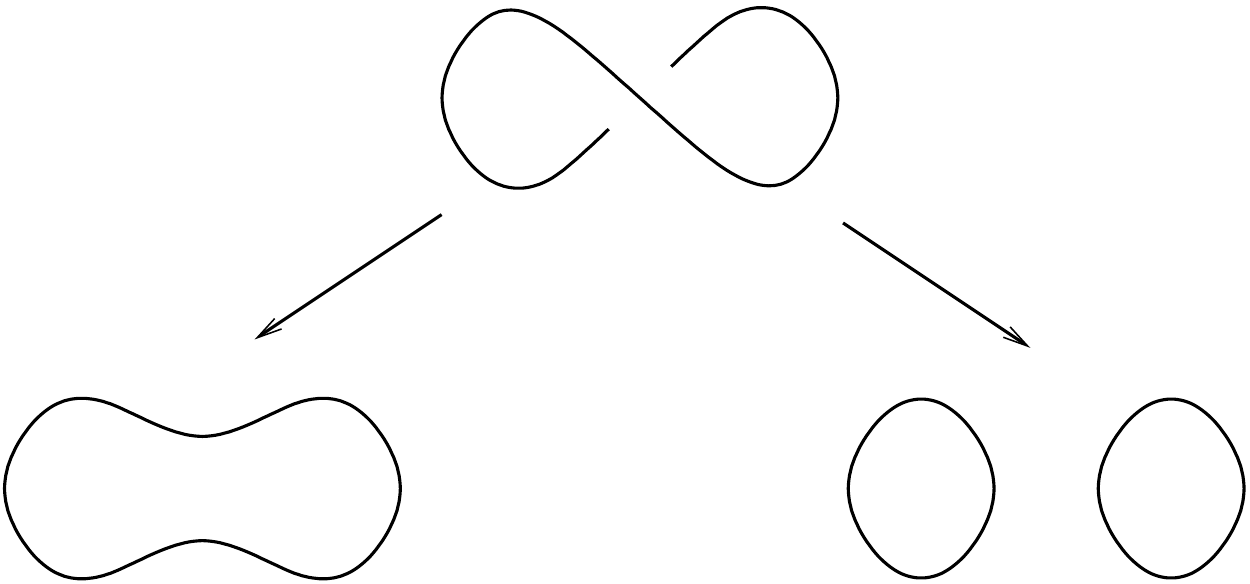}
\end{figure}

By analogy with the previous case, we want to lift the formula
for $\langle D \rangle$ to a complex
$$ 0 \lra \ A  \ \stackrel{\Delta}{\lra} \ A \otimes A\{-1\} \  \lra \  0 $$
with the differential denoted by $\Delta$, since it resembles
comultiplication. The degrees for basis elements of
these vector spaces are listed in the table below.

\begin{table}
\begin{tabular}{c|c|c}
Degree  &  Basis of $A$ & Basis of $A\otimes A\{ -1\}$  \\
\hline \hline
 $1 $  &  $X$  & $X \otimes X$    \\
\hline
 $0$   &             &  \\
\hline
$-1$   & $\mathbf{1}$ &  $ X\otimes \mathbf{1}, \mathbf{1}\otimes X$  \\
\hline
 $-2$  &             &  \\
\hline
$-3$   &   & $ \mathbf{1} \otimes \mathbf{1} $
\end{tabular}
\end{table}

\newpage

The differential should preserve the internal grading.
This and other natural conditions, including
that the homology of this complex should be isomorphic to $A$,
modulo a bigrading shift, leads to a formula for $\Delta$:

\begin{table}[!htbp]
\begin{tabular}{lll}
$X $  & $\lra$ & $X \otimes X$  \\
$\mathbf{1}$   & $\lra$ & $\mathbf{1}\otimes X + X\otimes \mathbf{1} $
\end{tabular}
\end{table}

Maps $M$ and $\Delta$ turn $A$ into a commutative Frobenius algebra -
a commutative unital algebra equipped with a nondegenerate symmetric
trace form. Such algebras are in a bijection with 2-dimensional
topological quantum field theories, that is, monoidal functors
from the category of 2-dimensional cobordisms between 1-manifolds
to the category of vector spaces. $A$ is, in addition, graded,
and the degree of the map associated to a cobordism $S$ equals
minus the Euler characteristic of $S$.

An arbitrary
diagram $D$ with $n$ crossing has $2^n$ resolutions into
crossingless diagrams. With each
resolution we associate $A^{\otimes k}$, where $k$ is the number of
circles in it. These powers of $A$ can be naturally placed into
the vertices of an $n$-dimensional cube. Every edge of a cube
corresponds to an elementary modification of a resolution,
converting a $k$-circle planar diagram to a $(k\pm 1)$-circle
diagram. This modification (which can be realized by a cobordisms
between the circle diagrams) induces a map between corresponding
tensor powers of $A$ (the map is either $M$ or $\Delta$ times
identity on the remaining circles).

Furthermore, for every square face of the cube, two compositions
of maps assigned to its edges commute, due to $M$ and $\Delta$
being structure maps of a two-dimensional TQFT.

Taking direct sums of tensor powers of $A$ (with the
suitably shifted internal grading) along the hyperplanes
orthogonal to the main diagonal of the cube  and defining the differential
to be a signed sum of the edge maps of the cube nets us
a complex of graded vector spaces $\overline{C}(D)$.
The Euler characteristic of $\overline{C}(D)$ is the Kauffman
bracket $\langle D \rangle$.

One then shifts the bigrading of $\overline{C}(D)$ to match
the coefficient in the formula (\ref{jones-kauffman})  and gets the complex $C(D)$
associated to a planar diagram $D$.
Its homology is denoted $H(D)$ and carries a bigrading.

\begin{theorem}
$H(D)$ depends only on the underlying link $L$ and not on its diagram $D$. The Euler characteristic of $H(D)$ is the Jones polynomial $J(L)$.
\end{theorem}

Denoting $H(D)$ by $H(L)$, we obtain a homology theory of links in
$\R^3$. It can be naturally thought of as a categorification
of the Jones polynomial (the term \emph{categorification} was
originally introduced by Louis Crane and Igor Frenkel in a
related context~\cite{CF}).

Since its discovery, this homology theory of links has been greatly developed and thoroughly understood by many people.
Below is a very brief list of some of the benefits and structures
 stemming from the homology theory $H$.

\begin{itemize}
\item $H$ contains large amount of information about knots and links and gives rise to
new structural relations between low-dimensional topology and algebra.
\item $H$ is functorial and extends to an invariant of link cobordisms (Jacobsson, Bar-Natan, Khovanov, Clark-Morrison-Walker).
\item $H$ carries homological operations (Lipshitz-Sarkar, Kriz-Kriz-Po).
\item Other link polynomials have also been categorified,
including the  Alexander polynomial (Ozsv\'ath-Rasmussen-Szab\'o), the HOMFLYPT polynomial and more general
Reshetikhin-Turaev invariants.
\item $H$ relates to several areas of math (geometric representation theory,
symplectic topology, algebraic geometry, the Langlands program).
\item $H$ appears in mathematical physics (Gukov-Schwarz-Vafa, Witten).
\end{itemize}

We would like to emphasize that the overall structure is built
in two steps. First step is constructing the Jones polynomial,
which we can think of as an example of linearization, in this
case going from links, topological objects (which admit a discrete
combinatorial interpretation) to elements of a linear structure,
the ring of Laurent polynomials. The second step is categorification,
lifting linear invariants of links to vector spaces and homology
groups. It is  possible that categorification can be
thought of as a kind of second linearization, but we will
not try to carefully phrase here what this might mean.

We conclude with an exercise and a problem for the reader.

\begin{exercise} Can the following structures and operations be interpreted as
linearizations?
\begin{itemize}
\item Passing from a topological space to the ring of continuous
functions on it.
\item Quantization.
\item Quantum computation.
\end{itemize}
\end{exercise}

\begin{problem} Discover new linearizations and categorifications
and develop them.
\end{problem}

\section{Acknowledgments}
While writing this paper, the author was partially supported
by the NSF grant DMS-1406065.
The paper grew out of the talks the author gave at the MAA MathFest in August 2015 in Washington DC and Knots in Washington in December 2015.


\vspace{0.2in} 
{\small 
\noindent 
khovanov@math.columbia.edu \\
Department of Mathematics \\
Columbia University \\  
New York, NY 10027
}


\begin{thebibliography}{99}

\bibitem{CF} L.~Crane and I.~B.~Frenkel, Four dimensional topological quantum field theory, Hopf categories, and the canonical bases,
\emph{J.Math.Phys.} {\bf 35}, (1994) 5136-5154.

\bibitem{J} V.~F.~R.~Jones, A polynomial invariant for knots via von Neumann algebras,    \emph{Bull. Amer. Math. Soc. (N.S.)}
  {\bf 12}, no.~1 (1985), 103-111.

\bibitem{Ka} L.~Kauffman, State models and the Jones polynomial, Topology {\bf 26}, no.~3 (1987), 395–407.

\bibitem{K} M.~Khovanov, A categorification of the Jones polynomial.
\emph{Duke Math. J.} {\bf 101}, no.~3 (2000), 359-426.

\bibitem{KT} M.~Khovanov, Y.~Tian, Work in progress.

\bibitem{PSMW} L.~Page, S.~Brin, R.~Motwani, T.~Winograd, The PageRank citation ranking:
bringing order to the web. Technical report. Stanford InfoLab (1999).

\bibitem{N} S.~Nakamoto, Bitcoin: A peer-to-peer electronic cash system, www.cryptovest.co.uk
(2008).

\end{thebibliography}
\end{document}